\documentclass{article}
\usepackage{amsmath,amssymb,amsthm}
\usepackage{graphics}

\newtheorem{theorem}{Theorem}
\newtheorem{lemma}{Lemma}

\newcommand{\conv}{\mathop{\rm conv}\nolimits}
\newcommand{\dist}{\mathop{\rm dist}\nolimits}

\newtheorem{cor}{Corollary}

\title {AN EXTENSION OF DELSARTE'S METHOD. THE KISSING PROBLEM IN THREE AND FOUR DIMENSIONS}

\author {Oleg R. Musin \thanks{Institute for Math. Study of Complex Systems, Moscow State University, Moscow, Russia omusin@mail.ru}}

\begin{document}
\date{}
\maketitle




\section {Introduction}

The {\it kissing number} $k(n)$ is the highest number of equal nonoverlapping spheres in ${\bf R}^n$ that can touch another sphere of the same size. In three dimensions the kissing number problem is asking how many white billiard balls can 
{\em kiss} (touch) a black ball. 

The most symmetrical configuration, 12 billiard balls around another, is if the 12 balls are placed at positions corresponding to the vertices of a regular icosahedron concentric with the central ball. However, these 12 outer balls do not kiss each other and may all moved freely. So perhaps if you moved all of them to one side a 13th ball would possibly fit in?       

This problem was the subject of a famous discussion between Isaac Newton 
and David Gregory in 1694. (May 4, 1694; see interesting article \cite{Sz} for details of this discussion.) It is commonly said that Newton believed the answer was 12 balls, while Gregory 
thought that 13 might be possible. However, Bill Casselman \cite{Cas} found some puzzling features in this story.
 
This problem is often called the {\it thirteen spheres problem}. R. Hoppe \cite{Hop} thought he had solved the problem in 1874. But, Thomas Hales \cite{Hales}  in 1994 published analysis of
Hoppe's mistake (see also \cite{Sz1}). 
Finally this problem was solved by Sch\"utte and van der Waerden in 1953 \cite{SvdW2}. A subsequent two-pages sketch of an elegant proof was given  by Leech \cite{Lee} in 1956. No much doubts that Leech's proof is correct, but there are gaps in his exposition, many involved sophisticated spherical trigonometry. (Leech's proof was presented in the first edition  of the well known book by Aigner $\&$ Ziegler \cite{AZ}, the authors removed this chapter from the second edition because 
a complete proof to include so much spherical trigonometry.)
The thirteen spheres problem continues to be of interest, new proofs have been published in the last few years by Wu-Yi Hsiang \cite{Hs}, K\'aroly B\"or\"oczky \cite{Bor}, and Kurt Anstreicher \cite{Ans}.

Note that $k(4)\geqslant 24$.
Indeed, the unit sphere in ${\bf R}^4$ centered at $(0,0,0,0)$ has 24 unit spheres around it, centered at the points $(\pm\sqrt{2},\pm\sqrt{2},0,0)$, with any choice of signs and any ordering of the coordinates. The convex hull of these 24 points yields a famous 4-dimensional regular polytope - the ``24-cell". Its facets are 24 regular octahedra. 

Coxeter proposed upper bounds on $k(n)$ in 1963  \cite{Cox}; for $n=4, 5, 6, 7,$ and 8 these bounds were 26, 48, 85, 146, and 244, respectively. Coxeter's bounds are based on the conjecture that equal size spherical caps on a sphere ${\bf S}^k$ can be packed no denser than $k+1$ spherical caps on ${\bf S}^k$ that simultaneously touch one another.  
B\"or\"oczky proved this conjecture in 1978 \cite{Bor1}.

The main progress in the kissing number problem in high dimensions was in the end of 1970's. Vladimir Levenshtein \cite{Lev2}, and independently Andrew Odlyzko and Neil Sloane \cite{OdS}, 
\cite[Chap.13]{CS} using Delsarte's method in 1979 proved that
$k(8)=240$, and $k(24)=196560$. This proof is surprisingly short, clean, and technically easier 
than all proofs in three dimensions. 

However, $n=8, 24$ are the only dimensions in which this method gives a precise result. For other dimensions (for instance, $n=3, 4$) the upper bounds exceed the lower. 
In  \cite{OdS} the Delsarte method was applied in dimensions up to 24 (see \cite[Table 1.5]{CS}). For comparison with the values of Coxeter's bounds on $k(n)$ for $n=4, 5, 6, 7,$ and 8 this method gives 25, 46, 82, 140, and 240, respectively. (For $n=3$  Coxeter's and Delsarte's  methods only gave $k(3)\leqslant 13$ \cite{Cox,OdS}.)
Kabatiansky and Levenshtein have found an asymptotic upper bound $2^{0.401n(1+o(1))}$ for $k(n)$ in 1978 \cite{Kab}. The lower bound $2^{0.2075n(1+o(1))}\; $ was found in \cite{Wyn}.

Improvements in the upper bounds on kissing numbers (for $n<24$) were rather weak during next years 
(\cite[Preface to Third Edition]{CS} gives a brief review and references). Arestov and Babenko \cite{AB1} proved that 
the bound $\; k(4)\leqslant 25\; $ cannot be improved using Delsarte's method. 
Hsiang \cite{Hs1} claims a proof of $k(4)=24.$ His work has not received yet a positive peer review.

If $M$ unit spheres kiss the unit sphere in ${\bf R}^n$, then the set of kissing points 
is an arrangement on the central sphere such that the (Euclidean) distance between any two points is at least 1. So the kissing number problem can be stated in other way: How many points can be placed on the surface of ${\bf S}^{n-1}$ so that the angular separation between any two points is at least $60^\circ$? 

This leads to an important generalization: a finite subset $X$ of ${\bf S}^{n-1}$ is called a {\it spherical $z$-code} if for 
every pair  $(x,y)$ of $X$ the scalar product $x\cdot y\le z.$ Spherical codes have many applications. The main application outside mathematics is in the design of signals for data transmission and storage. There are interesting applications to the numerical evaluation of $n$-dimensional integrals \cite[Chap.3]{CS}.

The Delsarte method (also known in coding theory as Delsarte's linear programming method, Delsarte's scheme, polynomial method) is described in \cite{CS,Kab}. 
Let $f(t)$ be a real polynomial such that $f(t)\le 0$ for $t\in [-1,z]$, 
the coefficients $c_k$'s in the expansion of $f(t)$ in terms of Gegenbauer polynomials $G_k^{(n)}$ are nonnegative, and $c_0=1.$ Then the maximal number of points in a spherical $z$-code in ${\bf S}^{n-1}$  is bounded by $f(1)$.
Suitable coefficients $c_k$'s can be found by the linear programming method \cite[Chapters 9,13]{CS}.

We found an extension of the Delsarte method in 2003 \cite{Mus}(see details in \cite{Mus2}), that allowed to prove the bound $k(4)<25$, i.e. $k(4)=24$. This extension yields also a proof $\; k(3)<13.$

The first version of these proofs was relatively short, but used a numerical solution of some nonconvex optimization problems. Later on \cite{Mus2} these calculations have been reduced to calculations of roots of polynomials in one variable. (This is not a big problem now, all computer algebra systems such as Maple, Mathematica, and Matlab can find roots. Also these calculations can be independently verified. If you have approximate values all roots of a polynomial, then you can check the existence of these roots by simple computations.) 

We present in this paper a new proof of the Newton-Gregory problem, an extension of Delsarte's method, and a proof that $k(4)=24.$ 

\section {The thirteen spheres problem: a new proof}

Let us recall the definition of {\em Legendre polynomials} $P_k(t)$ by recurrence formula: 
%
$$P_0=1,\;\; P_1=t,\;  P_2=\frac{3}{2}t^2-\frac{1}{2}, \ldots,\; P_k=\frac {2k-1}{k}\,t\,P_{k-1}-\frac{k-1}{k}\,P_{k-2};$$
$$  \mbox{ or equivalently \quad} \,P_k(t)=\frac{1}{2^k\,k!}\,\frac{d^k}{dt^k}(t^2-1)^k \quad \mbox{ (Rodrigues' formula)}.$$

%

\begin {lemma}
Let $X = \{x_1, x_2,\ldots, x_n\}$ be any finite subset of the unit sphere ${\bf S}^2$ in ${\bf R}^3$.
By $\phi_{i,j}=\dist(x_i,x_j)$ we denote the spherical (angular) distance between  $x_i$ and  $x_j.$ Then
$$\sum\limits_{i=1}^n \sum\limits_{j=1}^n P_k(\cos(\phi_{i,j})) \geqslant 0.    $$
\end {lemma}

Let 
$$f(t) = \frac{2431}{80}t^9 - \frac{1287}{20}t^7 + \frac{18333}{400}t^5 + \frac{343}{40}t^4 - \frac{83}{10}t^3 - \frac{213}{100}t^2+\frac{t}{10} - \frac{1}{200}. $$

\medskip


\begin {lemma} Suppose $X = \{x_1, x_2,\ldots, x_n\} \subset {\bf S}^2$. Then
$$S(X)=\sum\limits_{i=1}^n \sum\limits_{j=1}^n f(\cos(\phi_{i,j})) \geqslant n^2.$$
\end {lemma}


\begin {lemma} Suppose $X = \{x_1, x_2,\ldots, x_n\}$ is a subset of ${\bf S}^2$ such that the angular separation 
$\phi_{i,j}$ between any two distinct points $x_i, x_j$ is at least $60^{\circ}$.
Then
$$S(X)=\sum\limits_{i=1}^n \sum\limits_{j=1}^n f(\cos(\phi_{i,j})) < 13n.$$
\end {lemma}

\begin{theorem} $k(3)=12.$
\end{theorem}

\begin{proof} Suppose $X$ is a kissing arrangement on ${\bf S}^2$ with $n=k(3)$. Then $X$ is satisfying the assumptions in Lemmas 2, 3. Therefore,
$n^2 \leqslant S(X) < 13n$. From this follows $n<13,$ i.e. $n\leqslant 12.$ From other side we have $k(3)\geqslant 12$, then $n=k(3)=12.$ 
\end{proof}

We need the only one fact from spherical trigonometry, namely the {\em law of cosines}: 
$$\cos{\phi} = \cos{\theta_1}\cos{\theta_2}+\sin{\theta_1}\sin{\theta_2}\cos\varphi,$$ 
where for spherical triangle $ABC$ the angular lengths of its sides are $\theta_1,\theta_2,\phi$ and the angle between $AB, AC$ is $\varphi$ (Fig. 1). If $\varphi=90^{\circ}$, then $\cos{\phi} = \cos{\theta_1}\cos{\theta_2}$ (spherical Pythagorean theorem).

\begin{center}
\begin{picture}(320,140)(-230,-70)
\thinlines
\put(-90,-54){\line(0,1){108}}
\put(90,-54){\line(0,1){108}}
\put(-90,-54){\line(1,0){180}}
\put(-90,54){\line(1,0){180}}
\put(-90,-40){\line(1,0){180}}

\thicklines
\qbezier (-90,36)(-88,38)(-82,18)
\qbezier (-82,18)(-80,12)(-78,5)
\qbezier (-78,5)(-76,-1)(-74,-7)
\qbezier (-74,-7)(-72,-13)(-70,-18)
\qbezier (-70,-18)(-66,-27) (-62,-33)

\qbezier (-62,-33)(-60,-35)(-58,-37)
\qbezier (-58,-37)(-56,-38)(-54,-39)
\qbezier (-54,-39)(-52,-40)(-50,-41)
\qbezier (-50,-41)(-46,-42)(-35,-42)
\qbezier (-35,-42)(-30,-41)(-24,-41)
\qbezier (-24,-41)(-3,-42)(10,-41) 
\qbezier (10,-41)(16,-42)(20,-43)
\qbezier (20,-43)(26,-44)(30,-45)
\qbezier (30,-45)(35,-46)(40,-48) 
\qbezier (40,-48)(45,-49)(50,-48)
\qbezier (50,-48)(55,-46)(60,-40)
\qbezier (60,-40)(66,-34)(70,-24)
\qbezier (70,-24)(75,-10)(80,6)
\qbezier (80,6)(86,24)(89,46)

\thinlines
\multiput (-80,-54)(10,0){17}%
{\line(0,1){2}}
\put(-97,-63){$-1$}
\put(-55,-63){$-0.5$}
\put(8,-63){$0$}
\put(53,-63){$0.5$}
\put(82,-63){$0.8$}
\put(-98,-42){$0$}
\put(-83,-80){Fig. 2. The graph of the function $f(t)$}

\thicklines
\put(-230,-54){\line(1,2){40}}
\put(-150,-54){\line(-1,2){40}}
\put(-230,-54){\line(1,0){80}}
\thinlines
\qbezier (-185,16)(-190,14)(-195,16)
\put(-194,8){$\varphi$}
\put(-193,30){$A$}
\put(-242,-56){$B$}
\put(-148,-56){$C$}
\put(-221,-14){$\theta_1$}
\put(-168,-14){$\theta_2$}
\put(-191,-50){$\phi$}

\put(-200,-80){Fig. 1}

\end{picture}
\end{center}

\medskip

\medskip

{\bf Proof of Lemma 1.}

\medskip

This lemma easily follows from Schoenberg's theorem \cite{Scho} for Gegenbauer polynomials. 
Note that $P_k=G_k^{(3)}$. For completeness we give a proof of Lemma 1 here. 
In this proof we are using original Schoenberg's proof that based on the addition theorem for Gegenbauer polynomials.\footnote{Pfender and Ziegler\cite{PZ} give a proof as a simple consequence of the addition theorem for spherical harmonics. This theorem is not so elementary. The addition theorem for Legendre polynomials can be proven by elementary algebraic calculations.} 

The addition theorem for Legendre polynomials was discovered by Laplace and Legendre in 1782-1785:
$$P_k(\cos{\theta_1}\cos{\theta_2}+\sin{\theta_1}\sin{\theta_2}\cos{\varphi}) = \sum\limits_{m=0}^k c_{m,k}\,P_k^m(\cos{\theta_1})P_k^m(\cos{\theta_2})\,\cos{m\varphi}$$
$$= P_k(\cos{\theta_1})P_k(\cos{\theta_2})+
2\sum\limits_{m=1}^k\,\frac{(k-m)!}{(k+m)!}\,P_k^m(\cos{\theta_1})P_k^m(\cos{\theta_2})\,\cos{m\varphi},$$
where
$$P_k^m(t) = (1-t^2)^{\frac{m}{2}}\,\frac{d^m}{dt^m}P_k(t).$$
(See details in \cite{Car, Erd}.)

\medskip

\medskip

\begin{proof}
Let $X=\{x_1, \ldots, x_n\} \subset {\bf S}^2$ and $x_i$ has spherical (polar) coordinates $(\theta_i,\varphi_i)$. Then from the law of cosines we have:
$$\cos{\phi_{i,j}}=\cos{\theta_i}\,\cos{\theta_j}+\sin{\theta_i}\sin{\theta_j}\cos{\varphi_{i,j}},\quad
\varphi_{i,j}=\varphi_i-\varphi_j,$$
which yields
$$\sum\limits_{i,j}P_k(\cos{\phi_{i,j}})=\sum\limits_{i,j}\sum\limits_{m=0}^kc_{m,k}P_k^m(\cos{\theta_i})P_k^m(\cos{\theta_j})\cos{m\varphi_{i,j}}$$
$$ = \sum\limits_mc_{m,k}\sum\limits_{i,j}u_{m,i}u_{m,j}\cos{m\varphi_{i,j}}, \quad 
u_{m,i}=P_k^m(\cos{\theta_i}).$$ 

Let us prove that for any real $u_1, \ldots, u_n$
$$\sum_{i,j}u_iu_j\cos{m\varphi_{i,j}}\geqslant 0.$$

Pick $n$ vectors $y_1, \ldots, y_n$  in ${\bf R}^2$ with coordinates 
$y_i=(\cos{m\varphi_i}, \sin{m\varphi_i})$. If $y=u_1y_1+\ldots+u_ny_n,$ then
$$<y,y> \, = \, ||y||^2 = \sum_{i,j}u_iu_j\cos{m\varphi_{i,j}}\geqslant 0.$$
This inequality and the inequalities $c_{m,k}>0$ complete our proof.
\end{proof}

\medskip

\medskip

\medskip

\medskip

{\bf Proof of Lemma 2.}

\begin{proof}
The expansion of $f$ in terms of $P_k$ is
$$f = \sum\limits_{k=0}^9 {c_kP_k} = P_0 + 1.6P_1 + 3.48P_2 + 1.65P_3 + 1.96P_4 + 0.1P_5 + 0.32P_9.$$
We have $c_0=1,\; c_k \geqslant 0,\; k=1,2,\ldots, 9.\; $ Using Lemma 1 we get
$$ S(X)=\sum\limits_{k=0}^9 c_k \sum\limits_{i=1}^n \sum\limits_{j=1}^n P_k(\cos(\phi_{i,j}))\geqslant  
\sum\limits_{i=1}^n\sum\limits_{j=1}^n c_0P_0 =   n^2.$$
\end{proof}

\medskip

\medskip

{\bf Proof of Lemma 3.}
\begin{proof} {\bf 1.} 
The polynomial $f(t)$ satisfies the following properties (see Fig.2):

\noindent $(i)\; f(t)$ is a monotone decreasing function on the interval $[-1,-t_0];$
 
\noindent $(ii)\; f(t) < 0\;\;$ for $\; t\in (-t_0,1/2]; \\$ 
where $\; f(-t_0)=0, \; t_0 \approx  0.5907$.

These properties hold because $f(t)$ has the only one root $-t_0$ on $[-1, 1/2]$, and there are no zeros of the derivative $f'(t)$ (8th degree polynomial) on $[-1,-t_0].$

$$\mbox{Let} \; S_i(X):=\sum\limits_{j=1}^n f(\cos(\phi_{i,j})), \; \mbox{then} \; S(X)=\sum\limits_{i=1}^n S_i(X). \; \,  \mbox{From this follows} $$
if $ \; S_i(X)<13 \; $ for $ \; i=1, 2, \ldots, n, \; $ then $\; S(X)<13n$. 

We obviously have $\phi_{i,i}=0$, so $f(\cos{\phi_{i,i}})=f(1)$. Note that our assumption on  $X$ ($\phi_{i,j} \geqslant 60^{\circ}, \,  i \ne j$) yields
$\cos{\phi_{i,j}} \leqslant 1/2.$ Therefore, $\cos{\phi_{i,j}}$ lies in the interval [-1,1/2].
Since $(ii)$, if $\cos{\phi_{i,j}} \in [-t_0,1/2]$, then $f(\cos{\phi_{i,j}}) \leqslant 0$. Let $J(i):=\{j:\cos{\phi_{i,j}}\in [-1,-t_0)\}$. We obtain  
$$S_i(X) \leqslant T_i(X):=f(1)+\sum\limits_{j \in J(i)} f(\cos{\phi_{i,j}}).  \eqno (1)$$

Let $\theta_0=\arccos{t_0}, \theta_0 \approx 53.794^\circ.$ Then $j\in J(i)$ iff 
$\phi_{i,j}>180^\circ-\theta_0$, i.e.
$\theta_j<\theta_0$, where $\theta_j=180^\circ-\phi_{i,j}.$ In other words all $x_{i,j},\; j\in J(i)$ lie inside the circle of center $e_0$ and radius $\theta_0$, where $e_0=-x_i$ is the antipodal point to $x_i$.

\medskip

{\bf 2.} Let us consider on ${\bf S}^2$ points $e_0,y_1,\ldots,y_m$  such that
$$\phi_{i,j}=\dist(y_i,y_j)\geqslant 60^\circ \mbox{ for all }\; i\neq j,\quad \dist(e_0,y_i)\leqslant \theta_0 \; \mbox{ for } \; 1\leqslant i\leqslant m. \eqno (2)$$

Denote by $\mu $ the highest value of $m$ such that the constraints in $(2)$ define a non-empty set of points $y_1,\ldots,y_m.$ 

Suppose $\; 0\leqslant m\leqslant\mu\;$ and $Y=\{y_1,\ldots, y_m\}$ satisfies $(2)$. Let 
$$H(Y)=H(y_1,\ldots,y_m):=f(1)+f(-\cos{\theta_1})+\ldots+f(-\cos{\theta_m}), \; \theta_i=\dist(e_0,y_i)$$ 
$$h_m:=\max\limits_Y{H(Y)},\quad
h_{max}:=\max{\{h_0,h_1,\ldots,h_\mu\}}.$$

It is clear that $\; T_i(X)\leqslant h_m$, where $m=|J(i)|$.  From $(1)$ it follows that $ S_i(X) \leqslant h_m.$ Thus, if we prove that $h_{max} < 13$, then we prove Lemma 3. 

\medskip

{\bf 3.} Now we prove that $\mu \leqslant 4.$

\noindent Suppose $Y=\{y_1,\ldots,y_m\} \subset {\bf S}^2$ satisfies $(2)$. If $e_0$ is the North pole and $y_i$ has polar coordinates $(\theta_i,\varphi_i)$, then from the law of cosines we have:
$$\cos{\phi_{i,j}} = \cos{\theta_i}\cos{\theta_j}+\sin{\theta_i}\sin{\theta_j}\cos(\varphi_i-\varphi_j).$$
From $(2)$ we have $\cos{\phi_{i,j}}\leqslant 1/2$, then
$$\cos(\varphi_i-\varphi_j)\leqslant \frac{1/2-\cos{\theta_i}\cos{\theta_j}}{\sin{\theta_i}\sin{\theta_j}}.\eqno (3)$$
$$\mbox{Let }\quad Q(\alpha)=\frac{1/2-\cos{\alpha}\cos{\beta}}{\sin{\alpha}\sin{\beta}}, \; \; \mbox{ then } \; \;
Q'(\alpha)=\frac{2\cos{\beta}-\cos{\alpha}}{2\sin^2{\alpha}\sin{\beta}}.$$ 
From this follows, if $\; 0<\alpha, \beta\leqslant \theta_0$, then $ \cos{\beta}> 1/2$ (because $\theta_0<60^\circ$); so then $Q'(\alpha)> 0,$ and $Q(\alpha)\leqslant Q(\theta_0).$ Therefore,
$$\frac{1/2-\cos{\theta_i}\cos{\theta_j}}{\sin{\theta_i}\sin{\theta_j}} \leqslant
\frac{1/2-\cos^2{\theta_0}}{\sin^2{\theta_0}}=\frac{1/2-t_0^2}{1-t_0^2}.$$

Combining this inequality and (3), we get $$\cos(\varphi_i-\varphi_j)\leqslant 
\frac{1/2-t_0^2}{1-t_0^2}.$$ 
Note that $\arccos((1/2-t_0^2)/(1-t_0^2)) \approx 76.582^\circ>72^\circ$. Then $m\leqslant 4$ because no more than four points can lie in an unit circle with the minimum angular separation between any two points greater than $72^\circ$.

\medskip

{\bf 4.} Now we have to prove that $h_{max}=\max{\{h_0,h_1,h_2,h_3,h_4\}} < 13.\\$
We obviously have $h_0=f(1)=10.11<13$.

From $(i)$ follows that $f(-\cos{\theta})$ is a monotone decreasing function in $\theta$ on 
$[0,\theta_0].$  
Then for $m=1:\; H(y_1)=f(1)+f(-\cos{\theta_1})$ attains its maximum at $\theta_1=0,$ 
$$h_1=f(1)+f(-1)=12.88<13$$ 

{\bf 5.}
Let us consider for $m=2,3,4$ an {\em optimal} arrangement
$\{e_0,y_1,\ldots,y_m\}$ in ${\bf S}^2$ that gives maximum of $H(Y)=h_m$.
Note that for optimal arrangement points $y_k$  cannot be shifted towards $e_0$ because in this case $H(Y)$ increases. 

For $m=2$ this yields: $e_0\in y_1y_2,$ and  $\dist(y_1,y_2)= 60^\circ.$ If  $e_0\notin y_1y_2,$
then whole arc $y_1y_2$ can be shifted to $e_0$. Also if $\dist(y_1,y_2)> 60^\circ,$ then
$y_1$ (and $y_2$) can be shifted to $e_0.$

For $m=3$ we prove that $\Delta_3=y_1y_2y_3$ is a spherical regular triangle with edge length $60^\circ$. As above, $e_0\in \Delta_3$, otherwise whole triangle can be shifted to $e_0$. Suppose $\dist(y_1,y_i)>60^\circ,\; i=2,3,$ then $\dist(y_1,e_0)$ can be decreased. From this follows that for any $y_i$ at least one of the distances $\{\dist(y_i,y_j)\}$ is equal to $60^\circ$. Therefore, at least two sides of 
$\Delta_3$ (say $y_1y_2$ and $y_1y_3$) have length $60^\circ$. Also $\dist(y_2,y_3)=60^\circ$, conversely  $y_3$ (or $y_2$, if $e_0 \in y_1y_3$) can be rotated about $y_1$ by a small angle towards $e_0$ (Fig.3).  

When $m=4$ first we prove that $\Delta_4=y_1y_2y_3y_4$ is a convex quadrangle. 
Conversely, we may assume that $y_4\in y_1y_2y_3$. 

The great circle that is orthogonal to the arc $e_0y_4$ divides ${\bf S}^2$ into two hemispheres: $H_1$ and $H_2$.  Suppose 
$e_0\in H_1$,   then at least one $y_i$ (say $y_3$) belongs $H_2$ (Fig.4).  
So the angle $\; \angle {e_0y_4y_3}\; $ greater  than $90^\circ$, then (again from the law of cosines)  
$\quad \dist(y_3,e_0)>\dist(y_3,y_4). \quad $ Thus,
$\\ \theta_3=\dist(y_3,e_0)>\dist(y_3,y_4)\geqslant 60^\circ>\theta_0 \quad - \; \, $ a contradiction.

Arguing as for $m=3$ it is easy to prove that $\Delta_4$ is a spherical equilateral quadrangle (rhomb) with edge length $60^\circ$. 

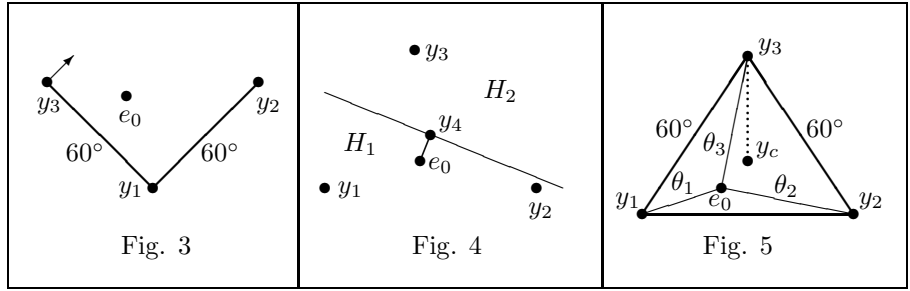
\begin{center}
\begin{picture}(320,110)(0,180)

\thicklines
\put(45,220){\line(-1,1){40}}
\put(45,220){\line(1,1){40}}

\thinlines

\put(5,260){\vector(1,1){10}}

%

\put(45,220){\circle*{4}}
\put(5,260){\circle*{4}}
\put(85,260){\circle*{4}}
\put(35,255){\circle*{4}}
\put(32,245){$e_0$}
\put(32,218){$y_1$}
\put(1,250){$y_3$}
\put(84,250){$y_2$}
\put(12,230){$60^\circ$}
\put(63,230){$60^\circ$}

\put(33,195){Fig. 3}

\thicklines
\put(230,210){\line(2,3){40}}
\put(310,210){\line(-2,3){40}}
\put(230,210){\line(1,0){80}}

\thinlines
\multiput(270,233.5)(0,2.5){15}%
{\circle*{1}}
\put(260,220){\line(1,5){10}}
\put(260,220){\line(-3,-1){30}}
\put(260,220){\line(5,-1){50}}

\put(270,230){\circle*{4}}
\put(230,210){\circle*{4}}
\put(310,210){\circle*{4}}
\put(270,270){\circle*{4}}
\put(260,220){\circle*{4}}

\put(220,214){$y_1$}
\put(312,214){$y_2$}
\put(274,272){$y_3$}
\put(273,233){$y_c$}
\put(255,213){$e_0$}
\put(241,219){$\theta_1$}
\put(280,218){$\theta_2$}
\put(253,234){$\theta_3$}
\put(235,239){$60^\circ$}
\put(292,239){$60^\circ$}
\put(253,195){Fig. 5}

\put(-10,290){\line(1,0){340}}
\put(-10,182){\line(1,0){340}}

\put(-10,182){\line(0,1){108}}
\put(330,182){\line(0,1){108}}

\put(100,182){\line(0,1){108}}
\put(215,182){\line(0,1){108}}


\put(150,240){\circle*{4}}
\put(110,220){\circle*{4}}
\put(190,220){\circle*{4}}
\put(146,230){\circle*{4}}

\put(115,218){$y_1$}
\put(187,210){$y_2$}
\put(148,270){$y_3$}
\put(144,272){\circle*{4}}

\put(153,243){$y_4$}
\put(149,227){$e_0$}

\put(110,256){\line(5,-2){90}}
\put(150,240){\line(-2,-5){4}}
\put(170,254){$H_2$}
\put(117,234){$H_1$}
\put(143,195){Fig. 4}




\end{picture}
\end{center}

{\bf 6.} Now we introduce the function $F_1(\psi),$\footnote{For given $\psi$, the value $F_1(\psi)$
can be find as the maximum of the  9th degree polynomial  $\Omega(s)=\widetilde F_1(\theta,\psi), \, s=\cos{(\theta-\psi/2)},$ on the interval $[\cos(\theta_0-\psi/2),1].$} 
where $\psi\in [60^\circ,2\theta_0]$:
$$F_1(\psi):=\max\limits_{\psi/2\leqslant \theta \leqslant \theta_0}
\{\widetilde F_1(\theta,\psi)\}, \quad 
\widetilde F_1(\theta,\psi)=f(-\cos{\theta})+f(-\cos(\psi-\theta)).$$
So if $\, \dist(y_i,y_j)=\psi,\;$ then
$$f(-\cos{\theta_i})+f(-\cos{\theta_j})\leqslant F_1(\psi).\eqno (4)$$ 

Therefore,

$$H(y_1,y_2)\leqslant h_2=f(1)+F_1(60^\circ)\approx 12.8749 < 13.$$

\medskip

{\bf 7.} When $m=4, \; \Delta_4$ is a spherical rhomb. Let 
$\; d_1=\dist(y_1,y_3)$, and $\\ d_2=\dist(y_2,y_4),\; $ then 
$\cos(d_1/2)\cos(d_2/2)=1/2\,$ (Pythagorean theorem, the diagonals of $\Delta_4$ are orthogonal). So if
$\rho(s):=2\arccos[1/(2\cos(s/2))],$ then $d_1=\rho(d_2), \; d_2=\rho(d_1),\; \rho(90^\circ)=90^\circ.$

 Suppose  $d_1\leqslant d_2$. Since $\theta_i\leqslant \theta_0, \; d_2\leqslant 2\theta_0,\; $  then   $\\ \rho(2\theta_0)\leqslant d_1 \leqslant 90^\circ\leqslant d_2\leqslant 2\theta_0.\;$

Now we consider two cases: $\\1) \; \rho(2\theta_0) \leqslant d_1<77^\circ, \; $ and $\; 2) \; 77^\circ\leqslant d_1\leqslant 90^\circ.\\$ 
$1)\; F_1(\psi)$ is a monotone decreasing function in $\psi.\;$ Then $(4)$ implies
$$f(-\cos{\theta_1})+f(-\cos{\theta_3})\leqslant F_1(\rho(2\theta_0)), \; f(-\cos{\theta_2})+f(-\cos{\theta_4}) < F_1(\rho(77^\circ)),$$
so then
$$H(Y) < f(1)+F_1(\rho(2\theta_0))+F_1(\rho(77^\circ))\approx 12.9171<13.$$
$2)$ In this case we have
$$H(Y) \leqslant f(1)+F_1(77^\circ)+F_1(90^\circ)\approx 12.9182<13.$$

Thus, $\; h_4 <13.$

\medskip

{\bf 8.} Our last step is to show that $h_3<13.$\footnote{More detailed analysis shows $h_3\approx 12.8721, \; h_4\approx 12.4849.$}

Since $\Delta_3$ is a regular triangle, $H(Y)=f(1)+f(-\cos{\theta_1})+f(-\cos{\theta_2})+f(-\cos{\theta_3})$ is a symmetric function in $\theta_i$, so we can consider only the case
$\\ \theta_1\leqslant \theta_2\leqslant \theta_3\leqslant \theta_0.\\ $ In this case
$\; R_0\leqslant \theta_3\leqslant\theta_0,\; $ where $\; R_0=\arccos{\sqrt{2/3}}\approx  35.2644^\circ .\; $  (Note that the circumradius of $\Delta_3\, $ equals $\, R_0$.) 

Let  $y_c$ is the center of $\Delta_3$. Denote by $u$ the angle $\angle {e_0y_3y_c}$. Then (see Fig.5)
$$\cos{\theta_1}=\cos{60^\circ}\cos{\theta_3}+\sin{60^\circ}\sin{\theta_3}\cos{(R_0-u)},$$
$$\cos{\theta_2}=\cos{60^\circ}\cos{\theta_3}+\sin{60^\circ}\sin{\theta_3}\cos{(R_0+u)},$$ 
where $\; \angle {y_1y_3y_c}=\angle {y_2y_3y_c}=R_0,\, \quad 0\leqslant u\leqslant u_0=\arccos(\cot{\theta_3}/\sqrt{3})-R_0\\$ (if $\, u=u_0, \, $ then $\, \theta_2=\theta_3$).

For fixed $\theta_3=\psi,\; H(y_1,y_2)$ becomes  the polynomial of degree 9 in $s=\cos{u}$. Denote by 
$F_2(\psi)$ the maximum of this polynomial on the interval $[\cos{u_0},1]$. 

Let
$$\{\psi_1,\ldots,\psi_6\}=\{R_0,\, 38^\circ,\, 41^\circ,\, 44^\circ,\, 48^\circ,\, \theta_0\}.\;$$
It's clear that $F_2(\psi)$  is a monotone increasing function in $\psi$ on  $[R_0,\theta_0].\; $ 
From other side, $f(-\cos{\psi})$ is a monotone decreasing function in $\psi$.
Therefore for  $\theta_3\in [\psi_i,\psi_{i+1}]$ we have
$$H(Y)=H(y_1,y_2)+f(-\cos{\theta_3}) < w_i:=F_2(\psi_{i+1})+f(-\cos{\psi_i}).$$ 
Since,
$$\; \{w_1,\ldots,w_5\}\approx\{12.9425, 12.9648, 12.9508, 12.9606, 12.9519\},$$ 
we get $\; h_3<\max\{w_i\}<13.$

Thus, $h_m <13$ for all $m$ as required.
\end{proof}

\section {Delsarte's method}

Let $X = \{x_1, x_2,\ldots, x_M\}$ be any finite subset of the unit sphere ${\bf S}^{n-1}
\subset{\bf R}^n,\\ {\bf S}^{n-1}=\{x: x\in {\bf R}^n,$ $x\cdot x=||x||^2=1\}.$ 
From here on we will speak of $x\in {\bf S}^{n-1}$ alternatively of points in ${\bf S}^{n-1}$ or of vectors in ${\bf R}^n.$

By $\phi_{ij}$ we denote the spherical (angular) distance between  $x_i,\,  x_j.$ It is clear that for any real numbers $u_1, u_2,\ldots, u_M$  the relation
$$ ||\sum u_ix_i||^2 = \sum\limits_{i,j} \cos{\phi_{ij}}u_iu_j \ge 0$$ holds, or equivalently
the Gram matrix $T(X)$ is  positive semidefinite, where \\ $T(X)=(t_{ij}),$ $t_{ij}=\cos{\phi_{ij}}=x_i\cdot x_j.$

Schoenberg \cite{Scho} extended this property to Gegenbauer (ultraspherical)  polynomials $G_k^{(n)}$ of $t_{ij}.$ He proved that if $g_{ij}=G_k^{(n)}(t_{ij}),$ then the matrix $(g_{ij})$ is positive semidefinite.  
Schoenberg proved also that the converse holds: if $f(t)$ is a real polynomial and for any finite $X\subset{\bf S}^{n-1}$ the matrix $(f(t_{ij}))$ is positive semidefinite, then $f$ is a sum of $G_k^{(n)}$ with nonnegative coefficients.

Let us recall the definition of Gegenbauer polynomials. Suppose $C_k^{(n)}(t)$ be the polynomials  defined by the expansion
$$(1-2rt+r^2)^{1-n/2} = \sum\limits_{k=0}\limits^{\infty}r^kC_k^{(n)}(t).$$ 
Then the polynomials $G_k^{(n)}(t) = C_k^{(n)}(t)/C_k^{(n)}(1)$ are called {\it Gegenbauer} or {\it ultraspherical} polynomials. (So the normalization of $G_k^{(n)}$ is determined by the condition $G_k^{(n)}(1)=1.$)

Also the Gegenbauer polynomials $G_k^{(n)}$ can be defined by recurrence formula:
$$G_0^{(n)}=1,\;\; G_1^{(n)}=t,\; \ldots,\; G_k^{(n)}=\frac {(2k+n-4)\,t\,G_{k-1}^{(n)}-(k-1)\,G_{k-2}^{(n)}} {k+n-3}$$

They are orthogonal on the interval $[-1,1]$ with respect to the weight function $\rho(t)=(1-t^2)^{(n-3)/2}$ (see details in \cite{Car,CS,Erd,Scho}). In the case $n=3,\; G_k^{(n)}$ are Legendre polynomials $P_k,$ and $G_k^{(4)}$ are  Chebyshev polynomials of the second kind (but with a different normalization than usual, $U_k(1)=1$),
$$ G_k^{(4)}(t)=U_k(t) = \frac {\sin{((k+1)\phi)}}{(k+1)\sin{\phi}}, \quad t=\cos{\phi}, \quad k=0,1,2,\ldots$$

For instance, $\;\; U_0=1,\;\;\; U_1=t,\;\;\; U_2=(4t^2-1)/3,\;\;\; U_3=2t^3-t,\\ U_4=(16t^4-12t^2+1)/5,\;\ldots,\;  
U_9=(256t^9-512t^7+336t^5-80t^3+5t)/5.$

\medskip

Let us now prove the bound of Delsarte's method. If a matrix $(g_{ij})$ is positive semidefinite,  then for any real 
$u_i$ the inequality  $\sum{g_{ij}u_iu_j}\geqslant 0$ holds, and then for $u_i=1,$ we have $ \sum\limits_{i,j}{g_{ij}}\ge 0.$ 
Therefore, for $g_{ij}=G_k^{(n)}(t_{ij}),$ we obtain
$$\sum\limits_{i=1}^M\sum\limits_{j=1}^M {{G_k^{(n)}(t_{ij})}} \geqslant 0 \eqno (3.1)$$

Suppose  $$f(t)=c_0G_0^{(n)}(t)+\ldots +c_dG_d^{(n)}(t),\; \mbox{ where }\; c_0\geqslant 0,  \ldots,\, c_d\geqslant 0. \eqno (3.2) $$
Let  $S(X)=\sum\limits_i\sum\limits_j{f(t_{ij})}.$ Using $(3.1),$ we get
$$S(X)=
\sum\limits_{k=0}^d\sum\limits_{i=1}^M\sum\limits_{j=1}^M {c_kG_k^{(n)}(t_{ij})}\geqslant  
\sum\limits_{i=1}^M\sum\limits_{j=1}^M {c_0G_0^{(n)}(t_{ij})} =   c_0M^2. \eqno (3.3)$$

Let $X=\{x_0,\ldots,x_M\}\subset {\bf S}^{n-1}$ be a spherical $z$-code, i.e. for all $i\neq j,$ $t_{ij}=\cos{\phi_{ij}}=x_i\cdot x_j\leqslant z,$ i.e. $t_{ij}\in [-1,z]$ (but $t_{ii}=1$). 
Suppose $f(t)\leqslant 0 $ for $t\in [-1,z],\;$ then \; $S(X)=Mf(1)+2f(t_{12})+\ldots+2f(t_{M-1\,M}) \leqslant Mf(1).$  If we combine this with $(3.2),$ then for $c_0>0$ we get 
$$M \leqslant \frac {f(1)}{c_0} \eqno (3.4)$$

The inequality $(3.4)$ play a crucial role in the Delsarte method (see details in  \cite{AB1, AB2, Boyv, CS, Del1, Del2, Kab, Lev2, OdS}). If $z=1/2$ and $c_0=1$,  then $(3.4)$ implies $k(n)\leqslant f(1).$ In \cite{Lev2, OdS} Levenshtein, Odlyzko and Sloane  have found the polynomials $f(t)$ such that  $f(1)=240,$ when $n=8;$ and $f(1)=196\,560,$ when  $n=24.$ Then
$k(8)\leqslant 240, \; k(24) \leqslant 196\,560.$ When $n=8, 24$, there exist sphere packings  ($E_8$ and Leech lattices) with these kissing numbers. Thus $k(8)=240$ and $k(24)=196\,560.$  When $n=4,$ a polynomial $f$ of degree 9 with $f(1)=25.5585...$ was found in \cite{OdS}. This implies 
$24\leqslant k(4) \leqslant 25.$

\section {An extension of Delsarte's method.}

Let us now generalize the Delsarte bound $M\leqslant f(1)/c_0$.

\medskip

\noindent {\bf Definition.} Let $f(t)$ be any function on the interval $[-1,1]$. 
Consider on ${\bf S}^{n-1}$  points $y_0,y_1,\ldots,y_m$  
such that
$$y_i\cdot y_j\leqslant z \; \mbox{ for all }\; i\neq j,\quad f(y_0\cdot y_i)>0 \; \mbox{ for } \; 1\leqslant i\leqslant m. \eqno (4.1)$$

Denote by $\mu=\mu(n,z,f) $ the highest value of $m$ such that the constraints in $(4.1)$ define a non-empty set of points $(y_0,\ldots,y_m).$ 

Suppose $\; 0\leqslant m\leqslant\mu.\; $ Let 
$$H(Y)=H(y_0;y_1,\ldots,y_m):=f(1)+f(y_0\cdot y_1)+\ldots+f(y_0\cdot y_m),$$ 
$$h_m:=\max\limits_Y\{H(Y)\},\quad
h_{max}:=\max{\{h_0,h_1,\ldots,h_\mu\}}.$$

\noindent {\bf Remark.} $h_{max}$ depends on $n,\; z,$ and $f.$ 
Throughout  this paper it is clear what $f,\; n,$ and  $z$ are; so we denote by $h_{max}$ the value $h_{max}(n,z,f).$

\medskip

\begin{theorem} Suppose $X\subset{\bf S}^{n-1}\; $ is a spherical $z$-code, $\; |X|=M, \; $ and  
$\\f(t)=c_0G_0^{(n)}(t)+\ldots +c_dG_d^{(n)}(t),$ where  $\; c_0>0,\; c_1\geqslant 0,  \ldots,\, c_d\geqslant 0.$ Then 
$$ M \leqslant \frac {h_{max}}{c_0}=\frac{1}{c_0}\max\{h_0,h_1,\ldots,h_\mu\}.$$
\end{theorem}

\begin{proof} Since $f$ satisfies $(3.2)$, then $(3.3)$ yields
$$S(X)\geqslant c_0M^2.$$

Let $J(i):=\{j:f(x_i\cdot x_j)>0, \; j\neq i\},$ and $X(i)=\{x_j: j\in J(i)\}.$ 
Then
$$S_i(X)=\sum\limits_{j=1}^M {f(x_i\cdot x_j)}\leqslant f(1)+\sum\limits_{j\in J(i)} {f(x_i\cdot x_j)}=H(x_i;X(i))\leqslant h_{max},$$  
so then $$S(X)=\sum\limits_{i=1}\limits^M S_i(X)\leqslant Mh_{max}.$$

We have $c_0M^2\leqslant S(X) \leqslant Mh_{max},$ i.e. $c_0M\leqslant h_{max}\; $ as required.
\end{proof}

Note that $h_0=f(1).$ If $f(t)\leqslant 0$ for all $t\in [-1,z]$, then for a $z$-code $X$ we have 
$\mu=0,$ i.e. $h_{max}=h_0=f(1).$ Therefore, this theorem yields the Delsarte bound
$M\leqslant f(1)/c_0.$

\medskip

The problem of evaluating of $h_{max}$ in general case looks even more complicated than the upper bound problem for spherical $z$-codes. It is not clear how to find $\mu$? Here we consider this problem only for a very restrictive class of functions $\; f(t)$: 
$f(t)\leqslant 0 \; \mbox{ for } \; t\in [-t_0,z], \quad t_0>z\geqslant 0. $

Let us denote by $A(k,\omega)$ the maximal number of points in a spherical $s$-code $\Omega \subset {\bf S}^{k-1}$ of minimal angle $\; \omega, \; \cos{\omega}=s.$ (Note that $A(n,60^\circ)$ is the kissing number $k(n)$.)

\begin {theorem} Suppose $Y=\{y_1,\ldots,y_m\}$ is a spherical $z$-code in ${\bf S}^{n-1}$, and points $y_i$ lie inside the sphere of center $e_0$ and radius $\theta_0$, where $t_0=\cos{\theta_0}\geqslant z.$ 
Then
$$m\leqslant A\Bigl(n-1,\arccos{\frac{z-t_0^2}{1-t_0^2}}\Bigr).$$ 
\end {theorem}
\begin{proof} We have 
$\; \phi_{i,j}=\dist(y_i,y_j)\geqslant \delta=\arccos{z} \;\; \mbox{ for } \; i\neq j;\\$
$\theta_i=\arccos(e_0\cdot y_i)\leqslant \theta_0\; \; \mbox{ for } \; \; 1\leqslant i\leqslant m;\; $ and
  $\; \theta_0\leqslant \delta.$

Let $\Pi$ be the projection of $Y$ onto equator ${\bf S}^{n-2}$ from pole $e_0.$ Denote by  $\gamma_{i,j}$ the distances between points of $\Pi$ in ${\bf S}^{n-2}$.  Then from the law of cosines and the inequality $\cos{\phi_{i,j}}\leqslant z,$ we get
$$\cos{\gamma_{i,j}}=\frac{\cos{\phi_{i,j}}-\cos{\theta_i}\cos{\theta_j}}{\sin{\theta_i}\sin{\theta_j}}
\leqslant 
\frac{z-\cos{\theta_i}\cos{\theta_j}}{\sin{\theta_i}\sin{\theta_j}}
$$ 
$$\mbox{Let }\quad Q(\alpha)=\frac{z-\cos{\alpha}\cos{\beta}}{\sin{\alpha}\sin{\beta}}, \; \; \mbox{ then } \; \;
Q'(\alpha)=\frac{\cos{\beta}-z\cos{\alpha}}{\sin^2{\alpha}\sin{\beta}}.$$ 
From this follows, if $\; 0<\alpha, \beta\leqslant \theta_0$, then $ \cos{\beta}\geqslant z$ (because $\theta_0\leqslant \delta$); so then $Q'(\alpha)\geqslant 0,$ and $Q(\alpha)\leqslant Q(\theta_0).$ Therefore,
$$\cos{\gamma_{i,j}}\leqslant \frac{z-\cos{\theta_i}\cos{\theta_j}}{\sin{\theta_i}\sin{\theta_j}} \leqslant
\frac{z-\cos^2{\theta_0}}{\sin^2{\theta_0}}=\frac{z-t_0^2}{1-t_0^2}$$
that complete our proof.
\end{proof}

\begin{cor} Suppose $f(t)\leqslant 0 \; \mbox{ for } \; t\in [-t_0,z], \quad t_0\geqslant z\geqslant 0,$ then
$$\mu(n,z,f)\leqslant A\Bigl(n-1,\arccos{\frac{z-t_0^2}{1-t_0^2}}\Bigr).$$ 
\end{cor}
\begin{proof} The assumption on $\; f\; $ yields $\; f(y_0\cdot y_i)>0\; $ only if 
$$\theta_i=\dist(e_0,y_i)<\theta_0=\arccos{t_0},$$ where $e_0=-y_0$ is the antipodal point to $y_0.$
Therefore, this set of points $\{e_0,y_1,\ldots,y_m\}$ satisfies the assumptions in Theorem 3.
\end{proof}

The next claim will be applied to prove that $k(4)=24.$

\begin{cor} Suppose $f(t)\leqslant 0 \; \mbox{ for } \; t\in [-t_0,1/2], \quad t_0\geqslant 0.6058,\; $
then $\\ \mu=\mu(4,1/2,f)\leqslant 6.$
\end{cor}
\begin{proof} Note that for $\; t_0\geqslant 0.6058, \; \arccos[(1/2-t_0^2)/(1-t_0^2)]>77.87^\circ.$
So Corollary 1 implies $\mu(4,1/2,f)\leqslant A(3,77.87^\circ).$

Denote by $\varphi_k(M)$ the largest angular separation that can be attained in a spherical code on ${\bf S}^{k-1}$ containing $M$ points. In three dimensions the best codes and the values $\varphi_3(M)$ presently known for $M\leqslant 12$   and $M=24$ (see \cite{Dan,FeT,SvdW1}).
For instance, 
Sch\"utte and van der Waerden \cite{SvdW1} proved that 
$\varphi_3(5)=\varphi_3(6)=90^\circ$ and $\varphi_3(7)\approx 77.86954^\circ \; (\cos{\varphi_3(7)}=\cot{40^\circ}\cot{80^\circ})$. 

Since $\; 77.87^\circ>\varphi_3(7),\; $ then $\; A(3,77.87^\circ)<7,\; $ i.e. $\; \mu\leqslant 6.$
\end{proof}

Corollary 1 shows that if $t_0$ is close enough to 1, then $\mu$ is small enough. Then one gets relatively small -
dimensional optimization problems for computation of numbers $h_m$ for small $n$. If additionally $f(t)$ is a monotone decreasing function on $[-1,-t_0]$, then these problems can be reduced to low-dimensional optimization problems of a type that can be treated numerically. 

\section{Optimal sets for monotonic functions}
In this section we consider $f(t)$ that satisfies the monotonicity assumption:
$$f(t) \; \mbox{ is a monotone decreasing function on the interval }\; [-1,-t_0],$$ 
$$f(t)\leqslant 0 \; \mbox{ for } \; t\in [-t_0,z], \quad t_0>z\geqslant 0 \eqno (*) $$.

Consider on ${\bf S}^{n-1}$  points $y_0,y_1,\ldots,y_m$  that satisfy $(4.1)$.
Denote by  $\theta_k$ for $k>0$ the distance between $y_k$ and $e_0,$ where $e_0=-y_0$ is the antipodal point to $y_0.$ Then $y_0\cdot y_k=-\cos{\theta_k},$  and $H(Y)$ is represented in the form:
$$H(Y)=f(1)+f(-\cos{\theta_1})+\ldots+f(-\cos{\theta_m}). \eqno (5.1)$$

A subset $C$ of ${\bf S}^{n-1}$ is called (spherical)  {\it convex} if it contains, with every two nonantipodal points, the small arc of the great circle containing them. If, in addition, $C$ does not contain antipodal points, then $C$ is called strongly convex. The closure of a convex set is convex and is the intersection of closed hemispheres (see details in \cite{DGK}). If a subset $Z$ of ${\bf S}^{n-1}$ lies in a hemisphere, then the convex hull of $Z$ is well defined, and is the intersection of all convex sets containing $Z$. 

Suppose $f(t)$ satisfies $(*)$, then $Q_m=\{y_1,\ldots, y_m\}$ lies in the hemisphere 
of center $e_0.$ 
Denote by $\Delta_m$ the convex hull of $Q_m$ in ${\bf S}^{n-1},\;$ $\Delta_m=\conv{Q_m}.$ 

Now we consider an optimal arrangement of $Q_m$ for $H.\; $ Let $\; \delta=\arccos{z},\\ \phi_{i,j}=\dist(y_i,y_j), \; \tilde N(Q_m)=$ number of $\phi_{i,j}=\delta \;\; (y_i\cdot y_j=z).$

\medskip

\noindent{\bf Definition} We say that $Q_m$
is optimal if $H(Y)=h_m$. If optimal $Q_m$ is not unique up to isometry, then we call $Q_m$ as optimal if it has maximal $\tilde N(Q_m)$.

\medskip

The function $f(t)$ is  monotone decreasing on $[-1,-t_0].$ By $(5.1)$ it follows that the function $H(Y)$  increases whenever $\theta_k$ decreases. This means that for an optimal $Q_m$ no $y_k \in Q_m$ can be shifted towards $e_0.$

That yields 
$$e_0\in \Delta_m \eqno (5.2)$$
because in the converse case whole $Q_m$ can be shifted to $e_0.$

From this follows that for $m=1, \; e_0=y_1$. Thus
$$h_1=f(1)+f(-1).$$

It was proved in Section 2 that for $m=2:\; \dist(y_1,y_2)=\delta$, 
thus
$$h_2=f(1)+\max\limits_{\delta/2\leqslant \theta \leqslant \theta_0}
\{(f(-\cos{\theta})+f(-\cos(\delta-\theta))\}, \quad 
\theta_0=\arccos{t_0}.$$

It was also proved that $\Delta_3$ is a spherical regular triangle with edge length $\delta.$ 
Using similar arguments it's not hard to prove that for $n>3,\; \Delta_4\; $ is a spherical regular tetrahedra with edge length $\delta.$ \footnote{For $m\leqslant n,\; \Delta_m\; $ is a spherical regular simplex with edge length $\delta.$ 
In this paper we need just cases $m=3, 4.$}

Let $\Delta_m,\; m\leqslant n,$ is a spherical regular simplex with edge length $\delta,$ and  
$$\Omega_m = \{y:y\in\Delta_m, \; y\cdot y_k\geqslant t_0,\; 1\leqslant k\leqslant m\}.$$ Note that $\Omega_m$  is a convex set in ${\bf S}^{n-1}$. Let $$H_m(y)=f(1)+f(-y\cdot y_1)+\ldots+f(-y\cdot y_m).$$ Then $h_m$ is the maximum of $H_m(y)$
 on $ \Omega_m.$ 
$$h_m=\max\limits_{y\in \Lambda_m}{H_m(y)}, \quad 
 \Omega_m\subset\Delta_m\subset {\bf S}^{n-1}, \quad 2\leqslant m\leqslant\min(n,\mu). \eqno (5.3)$$

When $n>m$ any $y_k\in Q_m$ is a vertex of $\Delta_m.$ In other words, no $y_k$ that lies inside $\Delta_m.$ In fact, that has been proved in Section 2 (see {\bf 5}, Fig. 4). 

In the first version of the paper \cite{Mus2} has been claimed that $\\$
{\em for optimal $Q_m$ with $m > n$, for any $y_k\in Q_m$ there are at least $n-1$ distinct points in $Q_m$  at the distance of $\delta$ from $y_k.\\$}  However, Eiichi Bannai and Makoto Tagami found some gaps in our exposition. Most of them are related to ``degenerated" configurations. In this paper we need only the case $n=4, \; m=5$. For this case they verified each step of our proof, considered all  ``degenerated" configurations, and finally gave clean and detailed proof. I wish to thank Eiichi Bannai and Makoto Tagami for this work. Now this claim in general case can be considered only as conjecture.

\section{An algorithm for computation suitable polynomials $f(t)$}
In this section is presented an algorithm for computation  ``optimal" \footnote{Open problem: is it true that for given $t_0, d$ this algorithm defines $f$ with minimal $h_{max}$?}  polynomials $f$ such that 
$f(t)$ is a monotone decreasing function on the interval $[-1,-t_0],$ and  
$f(t)\leqslant 0 \; \mbox{ for } \; t\in [-t_0,z], \quad t_0>z\geqslant 0$. This algorithm based on our knowledge about optimal arrangement of points $y_i$ for given $m$. Coefficients $c_k$ can be found via discretization and linear programming; such method had been employed already by Odlyzko and Sloane \cite{OdS} for the same purpose.

Let us have a polynomial $f$ represented in the form $f(t)=1+\sum\limits_{k=1}\limits^d c_kG_k^{(n)}(t)$. We have the following constraints for $f$: (C1) $\;\; c_k\geqslant 0,\;\; 1\leqslant k\leqslant d$;\\ (C2) $\; f(a)>f(b)\;$ for $\; -1\leqslant a<b\leqslant -t_0$;\quad (C3) $\; f(t) \leqslant 0\;$ for $\; -t_0\leqslant t\leqslant z.$

When $m\leqslant n, \; h_m=\max{H_m(y)}$, $y\in \Lambda_m$. We do not know $y$ where $H_m$ attains its maximum, so for evaluation of $h_m$ let us use $y_c$ $-$ the center of $\Delta_m.$ All vertices $y_k$ of $\Delta_m$ are at the distance of $R_m$ from $y_c,$ where \\ $\cos{R_m}=\sqrt{(1+(m-1)z)/m}.$

When $m=2n-2, \; \Delta_m$ presumably is  a regular $(n-1)$-dimensional cross-polytope. (It is not proven yet.) In this case $\; \cos{R_m}=\sqrt{z}.$ 

Let $I_n=\{1,\ldots,n\}\bigcup \{2n-2\}, \;\; m\in I_n, \;\; b_m=-\cos{R_m},\; $ whence \\ $H_m(y_c)=f(1)+mf(b_m).\;\; $ If $F_0$ is such that $H_m(y) \leqslant E=F_0+f(1),\;$ then (C4) $\; f(b_m)\leqslant F_0/m, \;\; m\in I_n.\;$ 
A polynomial $f$ that satisfies (C1-C4) and gives the minimal $E$ (note that $E=F_0+1+c_1+\ldots+c_d=F_0+f(1)$ will become a lower estimate of $h_{max}$) can be found by the following 

\medskip

\centerline{\bf Algorithm.} 

\medskip

\noindent Input: $\; n,\; z,\; t_0,\; d,\; N.$

\noindent Output: $\; c_1,\ldots, c_d,\; F_0, \; E.$

\medskip

\noindent {\it First} replace (C2) and (C3) by a finite set of inequalities at the points\\ $a_j=-1+\epsilon j,\;\; 0\leqslant j \leqslant N, \;\; \epsilon=(1+z)/N:$ 

\medskip 

\noindent {\it Second} use linear programming to find $F_0, c_1,\ldots, c_d$ so as to minimize \\ 
$E-1=F_0+\sum\limits_{k=1}\limits^dc_k\;\;$ subject to the constraints
$$c_k\geqslant 0,\quad 1\leqslant k\leqslant d;\qquad \sum\limits_{k=1}\limits^dc_kG_k^{(n)}(a_j)\geqslant \sum\limits_{k=1}\limits^dc_kG_k^{(n)}(a_{j+1}), \quad a_j\in [-1,-t_0];$$
 $$1+\sum\limits_{k=1}\limits^dc_kG_k^{(n)}(a_j)\leqslant 0,\quad a_j \in [-t_0,z];\quad 
1+\sum\limits_{k=1}\limits^dc_kG_k^{(n)}(b_m)\leqslant F_0/m,\quad m\in I_n.$$ 

Let us note again that $E = \max\limits_{m\in I_n} {H_m (y_c)} \leqslant h_{max}$ here, and that $E = h_{max}$ only if $h_{max} = H_{m_0}(y_c)$ for some $m_0 \in I_n$.

\section {On calculations of $h_m$ for $m\leqslant n$}
Here we explain how to solve the optimization problem $(5.3)$.
Let $\Delta_m\subset {\bf S}^{m-1}$ is a spherical regular simplex with edge length $\delta=\arccos{z};\; y_i,\; i=1,\ldots,m,$ are the vertices of $\Delta_m; \; t_i=y\cdot y_i=\cos{\theta_i}\geqslant t_0=\cos{\theta_0};\; t_0>z;\; f(t)$ is a monotone decreasing function on the interval $[-1,-t_0];\; h_m$ is the maximum of $H_m(y)$ subject to the constraints $t_i\geqslant t_0;\; H_m(y)=f(1)+f(-y\cdot y_1)+\ldots+f(-y\cdot y_m).$ 

\medskip

\noindent{\em The first method.} 

\noindent $H_m(y)$ is a symmetric function in the variables $\theta_1, \ldots, \theta_m$. Then we can consider this problem only on the domain $\Lambda=\{y:\theta_m\leqslant\ldots\leqslant \theta_2\leqslant\theta_1\}.$ Note that $\Lambda$ is a spherical simplex. Let us consider a barycentric triangulation of this simplex such that the diameter of any simplex $\sigma_i$ of this triangulation is not exceed $\epsilon.$

It is easy to prove that for any $y_k, \; y\cdot y_k$ attains its maximum on $\sigma_i$ at some vertex of $\sigma_i.$ Denote this vertex by $y_{k,i}$. Let $I=\{i: y_{1,i}\cdot y_1\geqslant t_0\}.$ So for $i\in I$ we have
$$f(-y_{k,i}\cdot y_k)=\max\limits_{y\in \sigma_i} {\{f(-y\cdot y_k)\}},$$
then
$$h_m\leqslant \max\limits_{i\in I}{\Bigl\{\sum\limits_{k=1}^m{f(-y_{k,i}\cdot y_k)}\Bigr\}}.   $$

That yields a very simple method for calculation of $h_m$. For $f$ from Section 9 this method gives $h_3\approx 24.8345,\; h_4\approx 24.818.$

\medskip

\noindent{\em The second method.}

\noindent For $m\leqslant n$ the values $h_m$ can be calculated another way. We are using here  that 
$f(t)=f_0+f_1t+\ldots+f_dt^d$ is a polynomial. The first method is technically easier then the second one. However, the second method doesn't assume that $f$ is a monotone decreasing function on $[-1,-t_0]$, and it can be applied to functions without monotonicity assumption.

Let us consider $H_m(y)$ as the symmetric  polynomial $F_m(t_1,\ldots,t_m)$ in the variables 
$t_i=y\cdot y
_i: F_m(t_1,\ldots,t_m)=f(1)+f(-t_1)+\ldots+f(-t_m).$ Denote by $s_k=s_k(t_1,\ldots,t_m)$ the power sum $t_1^k+\ldots+t_m^k.$ Then $$F_m(t_1,\ldots,t_m)=\Psi_m(s_1,\ldots,s_d)=f(1)+mf_0-f_1\,s_1+\ldots+(-1)^df_d\,s_d.$$
  
From the fact that $\Delta_m$ is a spherical regular simplex follows
$$s_2=\sigma(s_1):=\frac{z}{(m-1)z+1}s_1^2+1-z. \eqno (7.1)$$

Any symmetric polynomial in $m$ variables can be expressed as a polynomial of $s_1,\ldots, s_m.$ Therefore, in the case $k>m$ the power sum $s_k$ is $R_k(s_1,\ldots,s_m).$ Combining this with $(7.1)$, we get $$\Psi_m(s_1,\sigma(s_1),s_3,\ldots,s_d)=\Phi_m(s_1,s_3,\ldots,s_m).$$ Therefore, we have 
$$h_m=\max{\Phi_m(s_1,s_3,\ldots,s_m)},\quad (s_1,s_3,\ldots,s_m)\in D_m\subset 
{\bf R}^{m-1}, $$ where $D_m$ is the domain in ${\bf R}^{m-1}$ defined by the constraints 
$t_i\geqslant t_0$ and $(7.1).$

Let us show now how to determine $D_m$ for $m>2.$ The equation $(7.1)$ defines the ellipsoid
$E: s_2=\sigma(s_1)$ in space $\{t_1,\ldots,t_m\}.$ Then $s_1=t_1+\ldots+t_m$ attains its maximum on $E$ at the point with $t_1=t_2=\ldots=t_m,$ and $s_1$ achieves its minimum on
$E\bigcap \{t_i\geqslant t_0\}$ at the point with $t_2=\ldots=t_m=t_0.$ From this follows $w_1\leqslant s_1\leqslant w_2,$ where
$$w_1=\frac{\sqrt{(p-t_0^2)\,(p-z^2)}+z\,t_0}{p} +(m-1)\,t_0,\quad p=\frac{1+(m-2)\,z}{m-1},$$  $$w_2=\sqrt{m\,(m-1)\,z+m}.$$

The equation $s_1=\omega$ gives the hyperplane, and the equation $s_2=\sigma(\omega)$ gives the $(m-1)$-sphere in space: $\{(t_1,\ldots,t_m)\}$. Denote by  $S(\omega)$ the $(m-2)$-sphere that is the intersection of these hyperplane and sphere. Let $l_k(\omega)$ be the minimum of $s_k$ on  
$S(\omega)\bigcap \{t_i\geqslant t_0\},$ and $v_k(\omega)$ is its maximum. Now we have
$$h_m=\max\limits_{s_1}\max\limits_{s_3}\ldots\max\limits_{s_m}
{{{\Phi_m(s_1,s_3,\ldots,s_m)}}}, \; \mbox{ where }$$
$$w_1\leqslant s_1\leqslant w_2,\quad l_k(s_1)\leqslant s_k\le v_k(s_1),\;\; k=3,\ldots,m. $$

For the polynomial $f$ from Section 9 (and Section 2) we can give more details about calculations of $h_m$ for $m=3,4.$

Let us consider the case $m=3$ with $d=9$. In this case $F_\omega(s_3)=\Phi_3(\omega,s_3)$ is a polynomial of degree 3 in the variable $s_3.$

\begin{lemma}Let $f$ be a 9th degree polynomial $f(t)=f_0+f_1t+\ldots+f_9t^9$ such that
$f_9>0,\; f_6=f_8=0,$ and 
$f_7>-15f_9/7.$ If $F_\omega'(s)\le 0$ at $s=l_3(\omega)$, then the function $F_\omega(s)$ achieves its maximum on the interval $[l_3(\omega),v_3(\omega)]$ at $s=l_3(\omega).$
\end{lemma}
\begin{proof} The expansion of $s_9$ in terms of $\; s_1^{i}s_2^{j}s_3^{k},\;\; i+2j+3k=9,\;$ is
$$s_9=\frac19s_3^3+s_3^2(\frac23s_1^3+s_2s_1)+
s_3(\frac38s_2^3-\frac38s_2^2s_1^2-\frac78s_2s_1^4+\frac{5}{24}s_1^6)+R(s_1,s_2).$$ The coefficient of $s_3^2s_1$ in $s_7$ equals $7/9$. 
Thus  $$F_\omega(s)=-s^3\,f_9/9-s^2\,(f_9\,\omega\,\sigma(\omega)+2f_9\,\omega^3/3-7f_7\,\omega/9)+
sR_1(\omega)+R_0(\omega).$$
 $F_\omega(s)$ is a cubic polynomial with negative coefficient of $s^3.$ Then $F_\omega(s)$ is  
a concave function for $s>r,$ where $r: F_\omega''(r)=0.$ Therefore, if $r<l_3(\omega)$, then
$F_\omega(s)$ is a concave function on the interval $[l_3(\omega),v_3(\omega)]$. 
$r<l_3(\omega)$ iff
$$B(\omega):=3l_3(\omega)+6\omega^3+9\omega\,\sigma(\omega)>-7\omega{f_7/f_9}.$$ This inequality holds for $t_0<-z\le0.$ Indeed, $$\omega\geqslant w_1\geqslant 1+2z,\quad \sigma\,(\omega)\ge 1, \quad l_3(\omega)>0;$$ so then $$B(\omega)>15\omega>-7\omega{f_7/f_9}.$$
 The inequality $F'_\omega(l_3(\omega))\leqslant 0$ implies that $F_\omega(s)$ is a decreasing function on the interval $[l_3(\omega),v_3(\omega)]$. 
\end{proof}

The polynomial $f$ from Section 9  satisfies the  assumptions in this lemma. Then
$\Phi_3(\omega,s)$ attains its maximum at the point  $s=l_3(\omega),$ i.e. at the point with
$t_1=t_2\geqslant t_3,$ or with $t_1\geqslant t_2\ge t_3=t_0.$ If $t_1=t_2\geqslant t_3,\;$ then $p(\omega)=\Phi_3(\omega,l_3(\omega))$ is a polynomial in $\omega.$ This polynomial is a decreasing function in the variable $\omega$ on the interval $t_3\geqslant t_0.$ Therefore, $p(\omega)$ achieves its maximum on this interval at the point with $t_3=t_0.$ 
The calculations show that for $f$ from Section 9 $h_3=\max{p(\omega)}\approx 24.8345,$ when $\theta_3=\theta_0,\;\; \theta_1=\theta_2\approx 30.0715^\circ.$ 
\begin{cor} Let f be the polynomial from Section 9, then $h_3\approx 24.8345.$
\end{cor}

Consider the function $F_\omega(s_3,s_4)=\Phi_4(\omega,s_3,s_4)$ on $S(\omega)$. Let $q_i \in S(\omega)$ and $\; q_1:$ $t_1=t_2>t_3=t_4,\;\; q_2: t_1=t_2=t_3>t_4,\; $ and  $\; q_3: t_1>t_2=t_3=t_4.$
\begin{lemma} 
Let $f$ be a 9th degree polynomial $f(t)=\sum{f_it^i}.$ If $f_9>0$ and $f_6=f_8=0,$ then the function $F_\omega(s_3,s_4)$ achieves its maximum on $S(\omega)$ with $\omega>1$ at one of the points $(s_3(q_i),s_4(q_i)), \; i=1,2,3.$ 
\end{lemma}

\begin{proof}
The expansion of $s_9$ in terms of $s_1^{i}s_2^{j}s_3^{k}s_4^{l}\;$ is
$$s_9=\frac{9}{16}s_4^2s_1+\frac{1}{9}s_3^3-\frac{1}{3}s_3^2s_1^3+\frac{3}{4}s_4s_3s_1+\frac{3}{8}s_4s_2s_1^3 -\frac{3}{8}s_3s_2^2s_1^2-\frac{1}{24}s_3s_1^6+R(s_1,s_2). $$ 
The coefficient of $s_3^2s_1$ in $s_7$ equals $0$. We have $f_6=f_8=0$, then \\ $F_\omega(s_3,s_4)=-f_9\,s_9+\ldots=-f_9(s_3^3/9-s_3^2\,\omega^3/3)+\ldots $ Therefore,
$$F_{33}=\frac{\partial^2F_\omega(s_3,s_4)}{\partial^2s_3}=-f_9(\frac23s_3-\frac{2}{3}\omega^3)=\frac{2f_9}{3}(\omega^3-s_3). $$
If $F_\omega(s_3,s_4)$ has its maximum on $S(\omega)$ at the point $x,$ and $x$ is not a critical point of $s_3$ on $S(\omega),$ then $F_{33}\le 0.$ From other side,
for all $t_i\in [0,1]$ and $s_1=\omega>1$ we have $s_3\leqslant \omega < \omega^3$, so then $F_{33}>0$. 
The function $s_3$ on $S(\omega)$  (up to permutation of labels) has critical points at $q_i, \; i=1,2,3.$
\end{proof}
\begin{cor} Let f be the polynomial from Section 9, then $h_4\approx 24.818.$
\end{cor}
\begin{proof} By direct calculations it can be shown that \\
$\; F_\omega(s_3(q_1),s_4(q_1))>F_\omega(s_3(q_i),s_4(q_i))$ for $i=2,3.\; $
 Then Lemma 5 implies $h_4=\max{p(\omega)},$ where 
$p(\omega)= F_\omega(s_3(q_1),s_4(q_1))=\Phi_4(\omega,s_3(q_1),s_4(q_1)).$

The polynomial $p(\omega)$ attains its maximum $\; h_4\approx 24.818\;$ at the point with $\; \theta_1=\theta_2\approx 30.2310^\circ,\;\; \theta_3=\theta_4\approx 51.6765^\circ.$
\end{proof}


\section {On calculations of $h_5$ in four dimensions}
Let us consider the case $n=4,\; m=5.$ For simplicity here we consider only the case $z=1/2$. Then $\delta=60^\circ$ and $\theta_0=\arccos{t_0}<60^\circ.$

Denote by $\Gamma_5$   the graph of the edges of $\Delta_5$ with length $60^\circ$ , where $Q_5$ is an optimal set. 
The degree of any vertex of $\Gamma_5$ is not less than 3 (see Section 5). This implies that at least one vertex of $\Gamma_5$ has degree 4. Indeed, if all vertices of $\Gamma_5$ are of degree 3, then the sum of the degrees equals 15, i.e. is not an even number. 
There exists only one type of $\Gamma_5$ with these conditions (Fig. 6). 

\begin{center}
\begin{picture}(320,140)(-80,-70)
\put(65,-65){Fig. 6}

\put(10,-20){\circle*{5}}

\put(90,-40){\circle*{5}}
\put(70,20){\circle*{5}}

\put(130,40){\circle*{5}}

\put(70,60){\circle*{5}}

\thicklines
\put(10,-20){\line(4,-1){80}}
\put(10,-20){\line(3,2){60}}

\put(70,20){\line(3,1){60}}
\put(90,-40){\line(1,2){40}}

\put(70,60){\line(-3,-4){60}}
\put(70,60){\line(0,-1){40}}
\put(70,60){\line(1,-5){20}}
\put(70,60){\line(3,-1){60}}

\thinlines
\multiput(90,-40)(-1,3){20}%
{\circle*{1}}

\put(70,-14){$\alpha$}

\put(-5,-21){$y_5$}
\put(97,-41){$y_2$}
\put(57,24){$y_4$}
\put(134,44){$y_3$}
\put(73,65){$y_1$}
 
\end{picture}
\end{center}

For fixed $\dist(y_2,y_4)=\alpha, \; Q_5$ is uniquely defined up to isometry. Therefore, we have the 1-parametric family $\Delta_5(\alpha)$ on ${\bf S}^3.\;$ If $\dist(y_3,y_5)=\beta$, then
$$2\cos{\alpha}\cos{\beta}+\cos{\alpha}+\cos{\beta}=0 \eqno (8.1)$$
The equation $(8.1)$ defines the function $\beta=\lambda(\alpha).$ Then $\alpha=\lambda(\beta)$, 
 $\lambda(90^\circ)=90^\circ$. 

For all $i$ we have $\; \dist(y_i,e_0)\leqslant\theta_0,\; $ then $$\dist(y_i,y_j)\leqslant 
\dist(y_i,e_0)+\dist(y_j,e_0)\leqslant 2\theta_0.$$

Suppose $\alpha\leqslant\beta$, then $(8.1)$ and the inequality $\beta\leqslant 2\theta_0$ yield
$$\alpha_0\leqslant\alpha\leqslant 90^\circ\leqslant\beta\leqslant 2\theta_0, \quad \alpha_0:=\max\{60^\circ,\lambda(2\theta_0)\}.$$

Let  $$H_5(y,\alpha)=f(1)+f(-y\cdot y_1(\alpha))+\ldots+f(-y\cdot y_5(\alpha)).$$ Then
$$h_5=\max\limits_{y,\alpha}\{H_5(y,\alpha)\},\;\; y\in{\bf S}^3,\;\; y\cdot 
y_k(\alpha)\geqslant t_0,\;\; 1\leqslant k\leqslant 5, \;\; \alpha_0\leqslant \alpha\leqslant 90^\circ \eqno (8.2) $$

We have four-dimensional optimization problem  $(8.2)$. 
Our first approach for this problem was to apply numerical methods \cite{Mus}. For the polynomial $f$ from Section 9 this optimization problem  was solved numerically by using 
the Nelder-Mead simplex method: $H_5(y,\alpha)$ 
achieves its maximum $h_5\approx 24.6856$ at $\alpha=60^\circ$ and $\;y\;$ with 
$\; \theta_1\approx 42.1569^\circ,\; \theta_2=\theta_4\approx 32.3025^\circ, \; \theta_3=\theta_5=\theta_0.$ 
(The similar approach for the case $n=4,\; m=6$ gives the 3-parametric family $\Delta_6(\alpha,\beta,\gamma)$, and for $f$ from Section 9: $h_6\approx 22.5205.$)

Note that $(8.2)$ is a nonconvex constrained optimization problem. In this case, the Nelder-Mead simplex method and other local improvements methods cannot guarantee finding a global optimum. It's possible (using estimations of derivatives) to organize computational process in such way that it gives a global optimum. However, such kind solutions are very hard to verify and some mathematicians do not accept such kind proofs. Fortunately, an estimation of $h_5$ can be reduced to discrete optimization problems.

Let $\dist(y_1,y)=\psi$, and $\; \Phi_{i,j}(y,\psi)=f(-y\cdot y_i)+f(-y\cdot y_j).$ It's clear that for fixed $\psi,\; \Phi_{i,j}(y,\psi)$ attains its maximum at some point that lies in the great 2-sphere that contains $y_1, y_i, y_j.$ 
Now we introduce the function $F(\psi,\gamma)$.\footnote{$F(\psi,60^\circ)=F_2(\psi)-f(1)$ (see Section 2, {\bf 8}, Fig. 5).} 
Suppose $y_1y_iy_j$ is a spherical triangle in ${\bf S}^2$ with $\dist(y_1,y_i)=\dist(y_1,y_j)=60^\circ,\; \dist(y_i,y_j)=\gamma$, denote by $F(\psi,\gamma)$ 
the maximum of $\Phi_{i,j}(y,\psi)$ on ${\bf S}^2$ subject to the constraints $y\cdot y_k\geqslant t_0, \; k=i,j.\;$ Then 
$\; \Phi_{i,j}(y,\psi)\leqslant F(\psi,\gamma),$ so then $\; \Phi_{2,4}(y,\psi)\leqslant F(\psi,\alpha),$ $\; \Phi_{3,5}(y,\psi)\leqslant F(\psi,\beta).$ Thus
$$h_5\leqslant f(1)+f(-\cos{\psi})+F(\psi,\alpha)+F(\psi,\lambda(\alpha)).\eqno (8.3)$$

Let $\alpha_0<\alpha_1<\ldots<\alpha_k<\alpha_{k+1}=90^\circ.$
It's easy to see that $F(\psi,\gamma)$ is a monotone decreasing function in $\gamma.$ That implies
for $\; \alpha\in[\alpha_i,\alpha_{i+1}]:\\ F(\psi,\alpha)\leqslant F(\psi,\alpha_i), \;
F(\psi,\lambda(\alpha))\leqslant F(\psi,\lambda(\alpha_{i+1})).\;$ Therefore, from $(8.3)$ follows
$$h_5\leqslant   f(1)+f(-\cos{\psi})+\max\limits_{0\leqslant i\leqslant k} {\{F(\psi,\alpha_i)+
F(\psi,\lambda(\alpha_{i+1}))\}}.\eqno (8.4)$$

Note that $(8.4)$ to reduce the dimension of the optimization problem $(8.1)$ from 4 to 2. It is not too hard to solve this problem in general case. However, the polynomial $f$ from Section 9 satisfies an additional assumptions that allowed to find a weak bound on $h_5$ even more easier.

Let us briefly explain how to check the following assumptions for $f$:

\noindent  $1)\; \Phi_{i,j}(y,\psi)$ achieves its maximum at one of the ends of the arc $\omega(\psi,\gamma)$, where
$\omega(\psi,\gamma):=\{y: y\in{\bf S}^2,\; \dist(y_1,y)=\psi,\; y\cdot y_\ell\geqslant t_0, \; \ell=i,j\};$  

\noindent $2)\; F(\psi,\gamma)$ is a monotone increasing function in $\psi.$

For given $\gamma$ $(\gamma=\dist(y_i,y_j))$ and $\psi$ the function $\Phi_{i,j}(y,\psi)$ becomes a polynomial $p(s)$ of degree $d$ on $[s_0,1]$, where $\; s=\cos{u}, \;$   $u=\angle y_iy_1y_c$, and $y_c$ is the center of $y_1y_iy_j$ (see Section 2, {\bf 8}). Then $1)$ holds iff $p'(s)$ has no roots on $(s_0,1)$, either if $s: p'(s)=0$, then $p''(s)>0.$

Using $1)$ it's easy to check $2).$ For the polynomial $f(t)$ from Section 9 if $\gamma>62.41^\circ$, then $p(s)$ achieves its maximum at $s=s_0$ (i.e. $\dist(y_j,y)=\theta_0)$, so it's clear that $2)$ holds. From other side if $\gamma<69.34^\circ$, then the arc $\omega(\psi,\gamma)$ lies inside the triangle $y_1y_iy_j$, therefore $F(\psi,\gamma)$ increases whenever $\psi$ increases.

Note that $1)$ gives us the explicit expression for $F(\psi,\gamma)=\max(p(s_0),p(1))$. For fixed $\gamma$ and $\psi\leqslant\psi_\ell$ from $2)$ follows $F(\psi,\gamma)\leqslant F(\psi_\ell,\gamma)$.

Denote by $\psi_{L(i)},\; \psi_{U(i)}$ the lower and upper bounds on $\psi$ that defined by the constraints 
$\; \alpha\in [\alpha_i,\alpha_{i+1}], \; y\cdot y_q\geqslant t_0, \; q=1,\ldots,5.\\ $ 
Let $\psi_{L(i)}=\psi_{i,0}< \psi_{i,1}<\ldots<\psi_{i,\ell}<\psi_{i,\ell+1}=\psi_{U(i)}$. Recall that
$f(-\cos{\psi})$ is a monotone decreasing function in $\psi.$ Then 2) and $(8.4)$ yield
$$h_5\leqslant   f(1)+\max\limits_{0\leqslant i\leqslant k}{\max\limits_{0\leqslant j\leqslant \ell}
{\{R_{i,j}\}}}, \eqno (8.5)$$
where
$$R_{i,j}=f(-\cos{\psi_{i,j}})+ F(\psi_{i,j+1},\alpha_i)\}+F(\psi_{i,j+1},\lambda(\alpha_{i+1})).$$

It's very easy to apply this method. Here we need just to calculate the matrix $(R_{i,j})$ and the maximal value of its entries gives the bound on $h_5.$ For $f$ from Section 9 and $t_0\approx 0.60794,\; \theta_0=\arccos{t_0}\approx 52.5588^\circ,\; f(-t_0)=0$, this method gives the bound $h_5<24.8434.$\footnote{R achieves its maximum at $\alpha=60^\circ,\; \psi\approx 30.9344^\circ.$ Note that this bound exceeds the tight bound on $h_5\approx 24.6856$ given by numerical methods.}  

Now we show how to find an upper bound on $h_6.$ Let $\{e_0,y_1,\ldots,y_6\}\in {\bf S}^3,$
 $H(y_1,\ldots,y_6)=f(1)+f(-\cos{\theta_1})+\ldots+f(-\cos{\theta_6}),$ where $\theta_i=\dist(e_0,y_i)$.
Suppose $\theta_1\leqslant\theta_2\leqslant\ldots\leqslant\theta_6.$ Now we prove that $\theta_6\geqslant 45^\circ.$ That can be proven as Corollary 2 (Section 4). Conversely, all $\theta_i<45^\circ$. In this case $t_{0*}=\cos{\theta_6}>1/\sqrt{2}$, and $\omega=\arccos{[(1/2-t_{0*}^2)/(1-t_{0*}^2)]}>90^\circ.$
But if $u>90^\circ$, then $A(3,\omega)\leqslant 4$ (see  \cite{SvdW1, FeT}) - a contradiction. (In fact we proved that $\theta_5\geqslant 45^\circ$ also.)

Let us consider two cases: (i) $\; \theta_0\geqslant \theta_6>50^\circ\quad $ (ii) $\; 50^\circ\geqslant \theta_6\geqslant 45^\circ.$ 

\noindent (i) $H(y_1,\ldots,y_6)=H(y_1,\ldots,y_5)+f(-\cos{\theta_6}).\; $ We have
$$H(y_1,\ldots,y_5)\leqslant h_5<24.8434, \quad f(-\cos{\theta_6})<f(-\cos{50^\circ})\approx 0.0906,$$
then  $\; H(y_1,\ldots,y_6)<24.934.$

\noindent (ii) In this case all $\theta_i\leqslant 50^\circ.$ Therefore, we can apply $(8.5)$ for 
$\theta_0=50^\circ.$ This method gives $h_5(50^\circ)<23.9181,$ then $H(y_1,\ldots,y_5)\leqslant   
h_5(50^\circ)<23.9181,$ so then
 $$ H(y_1,\ldots,y_6)<23.9181+f(-\cos{45^\circ})\approx 23.9181+0.4533=24.3714$$

Thus $$h_6<\max\{24.934,24.3714\}=24.934$$

\medskip

\medskip

\section {$k(4)=24$}
For $n=4, \; z=\cos{60^\circ}=1/2$ we apply this extension of Delsarte's method with
$$f(t) = 53.76t^9 - 107.52t^7 + 70.56t^5 + 16.384t^4 - 9.832t^3 - 4.128t^2 -
 0.434t - 0.016$$
The expansion of $f$ in terms of $U_k=G_k^{(4)}$ is
$$f = U_0 + 2U_1 + 6.12U_2 + 3.484U_3 + 5.12U_4 + 1.05U_9$$
The polynomial $f$ has two roots on $[-1,1]$: $t_1=-t_0, \; t_0\approx 0.60794, \; t_2=1/2$,\\  
$f(t)\leqslant 0\;$ for $\;t\in [-t_0,1/2],$ and $f$ is a monotone decreasing function on the interval $[-1,-t_0].$ The last property holds because there are no zeros of the derivative $f'(t)$ on  $[-1,-t_0]$. Therefore, $f$ satisfies $(*)$ for $z=1/2.$ 

\medskip

\noindent {\bf Remark.} The polynomial $f$ was found by using the algorithm in Section 6.
This algorithm for $n=4,$ $z=1/2,$ $d=9,$ $N=2000,$ $t_0=0.6058$ gives $E\approx 24.7895.$ For the polynomial $f$ the coefficients $c_k$ were changed to ``better looking" ones with $E\approx 24.8644.$

\begin{center}
\begin{picture}(320,200)(-160,-110)
\thinlines
\put(-135,-80){\line(0,1){160}}
\put(135,-80){\line(0,1){160}}
\put(-135,-80){\line(1,0){270}}
\put(-135,80){\line(1,0){270}}
\put(-135,-60){\line(1,0){270}}

\thicklines
\qbezier (-135,54)(-132,51)(-129,45)
\qbezier (-129,45)(-126,37)(-123,27)
\qbezier (-123,27)(-120,18)(-117,8)
\qbezier (-117,8)(-114,-2)(-111,-11)
\qbezier (-111,-11)(-108,-19)(-105,-27)
\qbezier (-105,-27)(-99,-40) (-93,-49)
\qbezier (-93,-49)(-90,-52)(-87,-55)
\qbezier (-87,-55)(-84,-57)(-81,-59)
\qbezier (-81,-59)(-78,-60)(-75,-61)
\qbezier (-75,-61)(-69,-62)(-53,-62)
\qbezier (-53,-62)(-45,-61)(-37,-61)
\qbezier (-37,-61)(-10,-61)(15,-61) 
\qbezier (15,-61)(23,-62)(30,-63)
\qbezier (30,-63)(38,-65)(45,-67)
\qbezier (45,-67)(53,-69)(60,-71) 
\qbezier (60,-71)(68,-72)(75,-71) 
\qbezier (75,-71)(83,-68)(90,-60)
\qbezier (90,-60)(98,-50)(105,-35)
\qbezier (105,-35)(113,-16)(120,8)
\qbezier (120,8)(128,36)(134,69)

\thinlines
\multiput (-120,-80)(15,0){17}%
{\line(0,1){2}}
\multiput (-135,-40)(0,20){6}%
{\line(1,0){2}}
\put(-143,-90){$-1$}
\put(-119,-90){$-0.8$}
\put(-89,-90){$-0.6$}
\put(-59,-90){$-0.4$}
\put(-29,-90){$-0.2$}
\put(13,-90){$0$}
\put(40,-90){$0.2$}
\put(70,-90){$0.4$}
\put(100,-90){$0.6$}
\put(130,-90){$0.8$}
\put(-143,-62){$0$}
\put(-143,-42){$1$}
\put(-143,-22){$2$}
\put(-143,-2){$3$}
\put(-143,18){$4$}
\put(-143,38){$5$}
\put(-143,58){$6$}
\put(-150,-82){$-1$}
\put(-78,-110){Fig. 7. The graph of the function $f(t)$}

\end{picture}
\end{center}


%
\medskip

We have $t_0>0.6058.$ Then Corollary 2 gives $\mu\leqslant 6.$ 
Consider all $m\leqslant 6.$ 
$$h_0=f(1)=18.774,\quad h_1=f(1)+f(-1)=24.48.$$ 
$$h_2=f(1)+\max\limits_{30^\circ\leqslant \theta \leqslant \theta_0}
\{(f(-\cos{\theta})+f(-\cos(60^\circ-\theta))\}\approx 24.8644,$$
where $\theta_0=\arccos{t_0}\approx 52.5588^\circ.$

Note that $h_2$ can be calculated by the same method as in Section 2. Here $h_2= f(1)+2f(-\cos{30^\circ})$ also.


In Sections 7, 8 have been shown that
$$h_3\approx 24.8345, \quad h_4\approx 24.818, \quad h_5<24.8434, \quad h_6<24.934.$$

\begin{theorem}
$\quad k(4)=24$
\end{theorem}

\begin{proof} Let $X$ be a spherical $1/2$-code in ${\bf S}^3$ with $M=k(4)$ points.  
The polynomial $f$ is such that
$h_{max}<25,$ then combining this and Theorem 2, we get \\$k(4)\leqslant h_{max} < 25.$   Recall that $k(4)\geqslant 24.$ Consequently, $k(4)=24.$
\end{proof}

\section {Concluding remarks }
The algorithm in Section 6 can be applied to other dimensions and spherical $z$-codes. If $t_0=1,$ then the algorithm gives the Delsarte method. $E$ is an estimation of $h_{max}$  in this algorithm.

Direct application of the method developed in this paper, presumably could lead to some improvements in the upper bounds on kissing numbers in dimensions 9, 10, 16, 17, 18 given in \cite[Table 1.5]{CS}. (``Presumably" because the equality $\; h_{max}=E\; $ is not proven yet.)

In 9 and 10 dimensions Table 1.5 gives: $\\306\leqslant k(9)\leqslant 380,\quad 500\leqslant k(10)\leqslant 595.$\\
The algorithm gives:\\
$n=\;\,9:\; \deg{f}=11,\; E=h_1=366.7822,\; t_0=0.54;$\\
$n=10:\; \deg{f}=11,\; E=h_1=570.5240,\; t_0=0.586$.\\
For these dimensions there is a good chance to prove that $\\k(9)\leqslant 366,\; k(10)\leqslant 570.$

From the equality $k(3)=12$ follows $\varphi_3(13)<60^\circ.$ 
The method gives \\ $\varphi_3(13)<59.4^\circ$ ($\deg{f}=11$).
The lower bound on $\varphi_3(13)$ is $57.1367^\circ $ \cite{FeT}. Therefore, we have $57.1367^\circ\leqslant\varphi_3(13)<59.4^\circ.$

The method gives $\varphi_4(25)<59.81^\circ, \; \varphi_4(24) < 60.5^\circ.$ (This is theorem that can be proven by the same method as Theorem 4.) That improve the bounds:
$$\varphi_4(25)<60.79^\circ, \;\; \varphi_4(24) < 61.65^\circ \; \cite{Lev2} \; (\mbox{cf. } \cite{Boyv});
\;\; \varphi_4(24) < 61.47^\circ \; \cite{Boyv};$$
$$ \varphi_4(25)<60.5^\circ, \quad \varphi_4(24) < 61.41^\circ \; \cite{AB2}.$$
Now in these cases we have $$\quad 57.4988^\circ\ < \varphi_4(25) < 59.81^\circ,
\quad 60^\circ \leqslant \varphi_4(24) <  60.5^\circ.$$

\medskip

For all cases that were considered ($z\leqslant 0.6$) this method 
gives better bounds than Fejes T\'oth's bounds for $\varphi_3(M)$  \cite{FeT}  and Coxeter's bounds for all $\varphi_n(M)$ \cite{Cox}.
However, for $n=5,6,7$ direct use of this generalization of the Delsarte method does not give better upper bounds on $k(n)$ than the Delsarte method. It is an interesting problem to find better methods.

\medskip

\medskip

{\bf Acknowledgment.} I  wish to thank Eiichi Bannai, Ivan Dynnikov,
Dmitry Leshchiner, Sergei Ovchinnikov, Makoto Tagami and G\"unter Ziegler for helpful discussions and useful comments on this paper.

\end{document}